\documentclass[a4paper,12pt]{article}
\pdfoutput=1
\usepackage[paper=letterpaper,margin=1in]{geometry}
\usepackage{amsmath,amssymb,amsfonts,epsfig,url,color,todonotes}
\usepackage{titlesec}
\usepackage{amsmath,amsthm}

\renewenvironment{thebibliography}[1]{%
\begin{oldthebibliography}{#1}%
\setlength{\parskip}{0ex}%
\setlength{\itemsep}{0ex}%
}%
{%
\end{oldthebibliography}%
}

\titleformat{\section}[block]{\color{blue}\Large\bfseries\filcenter}{}{1em}{}
\titleformat{\subsection}[hang]{\bfseries}{}{1em}{}

\begin{document}

\begin{center}
{\Large \bf Erland Samuel Bring's \\``Transformation of Algebraic Equations''}
\end{center}
\medskip

\vspace{.4cm}
\centerline{
  {\large Alexander Chen}$^1$, 
  {\large Yang-Hui He}$^2$ \&
  {\large John McKay}$^3$
}
\vspace*{3.0ex}

\begin{center}
  {\it
    {\small
      ${}^{1}$
      Westminster School\\
      \qquad
          {\rm \url{alexander.chen@westminster.org.uk}}\\
    }
    \vspace*{1.5ex}
    {\small
      ${}^{2}$
      Merton College, University of Oxford, OX14JD, UK; \\
      Department of Mathematics, City, University of London, EC1V 0HB, UK; \\
      School of Physics, NanKai University, Tianjin, 300071, P.R.~China \\
        \qquad
            {\rm \url{hey@maths.ox.ac.uk}}\\
    }
    \vspace*{1.5ex}
    {\small
      ${}^{3}$ 
      Department of Mathematics and Statistics,\\
      Concordia University, 1455 de Maisonneuve Blvd.~West,\\
      Montreal, Quebec, H3G 1M8, Canada\\
      \qquad 
          {\rm \url{mckay@encs.concordia.ca}}\\
    }
  }
\end{center}

\vspace{1in}

\begin{abstract}
  We translate Erland Samuel Bring's treatise {\it Meletemata qu{\ae}dam Mathematica circa Transformationem {\AE}quationum Algebraicarum} (``Some Selected Mathematics on the Transformation of Algebraic Equations'') written as his Promotionschrift at the University of Lund in 1786, from its Latin into English, with modern mathematical notation.
  
  Bring (1736 - 98) made important contributions to algebraic equations and obtained the canonical form $x^5 + p x + q = 0$ for quintics before Jerrard, Ruffini and Abel. In due course, he realized the significance of the projective curve which now bears his name: the complete intersection of the homogeneous polynomials of degrees $1,2,3$ in $\mathbb{P}^4$, i.e.,
  $\left\{
  \sum\limits_{i=0}^4 z_i = \sum\limits_{i=0}^4 z_i^2 = \sum\limits_{i=0}^4 z_i^3 = 0
  \right\} \in \mathbb{C}[z_0,\ldots,z_4]$.
\end{abstract}

\newpage

\section{Introduction and Translators' Notes}
Within the illustrious history surrounding the Quintic Problem (q.~v.~\cite{harley,klein,ber,dm}), it is unavoidable that many tales should be eclipsed by the ultimate triumph and heart-breaking tragedy of Evariste Galois.
Dead more than a decade before Galois was born, Erland Samuel Bring (1736-98), late Professor of History at the University of Lund, certainly provided one such biography (cf.~e.g.,\cite{bio}).
Having read law, philosophy and history and subsequently practised in the first and lectured in the latter two, our polymath's true passion lay in mathematics.
Though with volumes in history and mathematics ranging from geometry, astronomy and commentaries on Euler, the magnum opus of Bring came in the form of his Promotionschrift in application for his professorship in mathematics -- when already professor of history -- at the University of Lund, then the Royal Academy of Queen Carolina, in 1786.

The great Sylvester begins \cite{sylvester} his treatise on Hamilton's numbers \cite{ham} with the words (as will be seen shortly, the dedication page of the thesis is by Somelius):
\begin{quote}
{\it
In the year 1786 Erland Samuel Bring, Professor at the University of Lund in Sweden, showed how by an extension of the method of Tschirnhausen it was possible to deprive the general algebraical equation of the 5th degree of three of its terms without solving an equation higher than the 3rd degree. By a well-understood, however singular, academical fiction, this discovery was ascribed by him to one of his own pupils, a certain Sven Gustaf Sommelius, and embodied in a thesis humbly submitted to himself for approval by that pupil as a preliminary to his obtaining his degree of Doctor of Philosophy in the University.
}
\end{quote}
It is thus remarkable that our humble professor of history has managed to bring the general quintic to the trinomial form
\begin{equation}\label{BJ}
x^5 + p x + q = 0 
\end{equation}
by solving a set of no higher than cubics whose analytic resolution are standard even for his time, using so-called Tschirnhaus transformations which date to the late 17th century.
The achievement, together with the thesis, was sadly much forgotten until later independently rediscovered by the Irish mathematician G.~Jerrard, {\it after} the works of Abel -- which Jerrard did not believe - showing the insolubility of the quintic by radicals.
The simple expression \eqref{BJ} is now called the Bring-Jerrard form of the quintic.
A history (including some parts of Bring's thesis in the original, p45, cit.~ibid.) was given by the Reverend R.~Harley \cite{harley} on this reduction.

It is amusing that Bring writes as the last sentence of \S III, for the general $m$-th degree polynomial that
\begin{quote}
{\it
\ldots in genere ad eliminandum terminum M:tum requiri ut resolvatur {\ae}quatio M:t{\ae} dignitatis, credo de reliquo quemcunque facile perspicere, unum tantummodo terminum ob eandem rationem, cujus supra mota mentio est, exterminari posse; adeo ut si ita sit, ut Algebraist{\ae} aliquando in spem venerint quamlibet {\ae}quationem ope transformationis a terminis suis intermedii liberatam resolvendi, fantendum omnino est, hujus exspectationi in hunc usque diem minime respondisse eventum nec unquam responsurum, quominus aliam transformandi methodum ac qu{\ae} hucusque obtinuit, olim adhibere contingat.
}
\end{quote}
(`` in general in order to remove the Mth term as required, an equation of degree M needs to be resolved. I believe, with regards to the remaining cases, anyone can easily discover that only one coefficient can be removed, by the same reasoning that has been mentioned already. If this is true, Algebrists finally will come to hope that any sort of equation, freed by the power of transformations, can be resolved by its coefficients, and it must be acknowledged wholly that another method of transformation other than this one has not obtained this to this day, but an expectation of this sort does not mean that an answer will never be found.'')
A forgiveable statement for the 18th century!

Yet Bring's optimism (beginning of \S.IV) remains sound advice to all generations of mathematicians:
\begin{quote}
{\it
Quamvis autem nihil fere ex iis, qu{\ae} ad Mathesin pertinent, intactum atque intentatum reliquerint amabilis hujus scienti{\ae} plusquam indefessi amantes atque cultores, ea tamen Algebr{\ae} particula, \ldots
}
\end{quote}
(However, no method belonging to the study of Mathematics should have been left behind, untried and extended, from this lovable science by its untiring lovers and worshippers, and this area of Algebra in particular, \ldots )

The wider significance of the Bring-Jerrard form was appreciated much later in the works of Hermite \cite{hermite} on elliptic functions and Klein's celebrated lectures \cite{klein} on the icosahedron (q.v.~an account especially in relation to Bring's curve in \cite{nash}; cf.~cite{BN,edge,weber}).
The key is that in reducing the general quintic
\begin{equation}
x^5 + a x^4 + b x^3 + c x^2 + d x + e = 0
\end{equation}
into the form of \eqref{BJ}, the five complex roots $z_{i=0,\ldots,4}$,  by Viet\`e, must be such that
\begin{equation}
0 = a = \sum\limits_{0 \leq i\leq 4} z_i \ , \quad
0 = b = \sum\limits_{0 \leq i<j \leq 4} z_i z_j \ , \quad
0 = c = \sum\limits_{0 \leq i<j<k \leq 4} z_i z_j z_k \ .
\end{equation}
Since symmetric polynomials in $z_i$ can be written interchangeably into the basis of elementary ones as above, or into sum of powers.
The above can be recast into
\begin{equation}\label{bring}
  \sum\limits_{i=0}^4 z_i = \sum\limits_{i=0}^4 z_i^2 = \sum\limits_{i=0}^4 z_i^3 = 0 \ .
\end{equation}
Klein's realization is that \eqref{bring} itself defines a complex algebraic curve, of genus 4, as a complete intersection in $\mathbb{P}^4$ with homogeneous coordinates $[z_0: \ldots : z_4]$.
This curve, which he dubbed {\bf Bring's curve} in honour of our protagonist, can be thought of as a moduli space of Bring-Jerrard quintics.
It has symmetry group $S_5$, the largest possible for a genus 4 curve.

More recently, the deep connection between Bring's curve to the icosahedron has been exploited in the context of dessins d'enfants \cite{ss},
sporadic groups and generalized Moonshine \cite{He:2015yoa}, as well as the $E_8$ Lie group and modular groups \cite{skip,yang}.
With this resurgence of interest it is therefore expedient to present, for the first time in other than the original Latin, the entirety of Bring's influential thesis, which hopefully is of interest to mathematicians and historians alike.

\subsection{Translaters' Synopsis of the Thesis}
\begin{description}
\item[I: ] A reminder of how the quadratic can be solved using linear substitutions;
\item[II: ] A similar treatment for the cubic;
\item[III: ] How one might attack the general equation of the $m$-th degree;
\item[IV: ] On real and complex condition of roots in solving the cubic;
\item[V: ] A general strategy for eliminating terms in the cubic;
\item[VI: ] Following the strategy of \S.~V, construct the specific elimination problem of removing the terms of degree 1 and 2;
\item[VII: ] Similar treatment as \S.~VI for the quartic;
\item[VIII: ] Further removal of the terms of  degrees 3 and 4 in the quartic;
\item[IX: ] Difficulties in removing three terms in the quartic;
\item[X: ] Setting up subsidiaries for the quintic;
\item[XI: ] Finally putting the quinitc into Bring-Jerrard form.
\end{description}

\subsection{Notes}
In translating the Latin, a few points may be of interest to the reader:
\begin{itemize}
\item Dignitatis:  The word for ``degree'' of the polynomial is ``dignitatis'', also, of course, meaning ``dignity'';
\item \&c, or ``etc'', is taken as ``$\ldots$'' in equations;
\item Greek is used only once in the thesis, \S III., p6, $\omega\zeta\ \varepsilon\nu \ \pi\alpha\rho o \delta\omega$, meaning ``in passing'' or ``just generally''.
\item Q.~E~.F.~: 
``quod erat faciendum'' (that which was to be done) used for the end of a {\it construction}, as opposed to the perhaps more familiar Q.~E.~D.~for ``demonstrandum'' used at the end of a {\it proof}; Indeed, one also has ``quad est absurdum'' - though not abbreviated - for proofs by contradiction;
\item Quartic equations are referred to as ``{\ae}quation{\ae} biquadratic{\ae}'';
\item ``S:{\ae} R:{\ae} M:tis'' is short for ``Serenissim{\ae} Regin{\ae} Majestatis'', used for the Swedish Queen (we thank Mr.~Julian Reid, archivist of Merton College, Oxford for this insight);
\item Terminus: The word for ``term'' in the polynomial is ``terminus'', it is also used to mean the coefficient of a term;
\item ``Variable'' corresponds to ``littera'', thus, e.g., ``littera $z$'' means ``the variable $z$'';
\item Various words such as ``real'', ``complex'', ``root'' are their obvious counterparts ``realis'', ``imaginaria'', ``radicus'';
\item When referring to variables, e.g., ``the value of $z$'', a Greek insertion is made as ``valor $\tau\omega\nu \ z$'' (most of the time
$\tau\omega\nu$ is shortened to $\tau8$) as  $\tau\omega\nu$ is Greek for ``of''.
\end{itemize}
Now, without further ado, we present, in English and modern notation, the thesis of Bring.

\newpage


\newpage

\section{{\LARGE Some Selected Mathematics on the Transformation of Algebraic Equations}}
which by the great consent of of the Faculty of Philosophy in the Royal Academy of Queen Carolina, Dr.~Erland Samuel Bring (Regius and Ordinary Professor of History) modestly presented by Sven Gustaf Sommelius (under the stipend of the King and of Palmcreutz of Lund) to the erudite examinations of the public, on the 15th day of December, 1786.

\section{Dedication}
To a man of great faith, of Her Most Serene Majesty, the noblest Lord Sven Lagerbring,
Doctor of Civil and Canon Law,
(most dignified counsellor of the Queen's Council, Celebrated Professor of History in the Queen Carolina Academy, most merited member of the Royal Academy, most venerated senior member of the Carolina Academy, my dearest master and uncle like an indulgent parent);

To a man of great faith, of Her Most Serene Majesty, the noblest Lord Charles A.~Hallenborg (Counsellor in the Council in the King's Dicastry of Bero, of great gravitas);

To a noble man and generous Lord Magnus Hallenborg (Celebrated accessory in the King's Dicastry of Junicopen, of great justice).

To a most noble and capable man Lord Adolf Hallenborg, the most noble Centurian of the Cavalry in her Royal Legion of South Scania;
To the Maecenians and most excellent patrons

I need neither to anxiously question the Maecenians nor appeal to doubting Patrons, under whose guardianship this Academic specimen was produced, as long as I and many others are allowed to pay our respects and adorations to you, most noble men, the most bold Maecenians and Patrons. 
Although it is truly most excellent, with regard to praises of you, to favour those more scholarly ones, with whose delights and pleasures you are truly considered; however this preference, shared by others as well as myself, I will honour with a silent loyalty, since there are currently other aspects which I wish to address more.

For I ought to remember publicly, lest I seem or be an ignorant child, the house of your people, in particular your forefathers, so many of whom were benefactors and so unusually magnanimous with favours, which are used by Parents/Founders for sermons/discussions. 
Therefore I put this on record, which is worth as much to me as a beloved memory, because it stimulates an emotion of loyalty to you, Maecenians and Great Patrons, and I decided to make testament to it with a deferential dedication. 
Therefore, I eagerly beseech and implore you to allow your most honourable name to be inscribed in this work, and more importantly, and again I beg you because I owe the fortunes of my beginning to your coming, that you allow this since it is ordained by destiny, and so it is also necessary for me to give back this account of my studies.

Not only do you allow me to offer these first fruits of labour but you have no issue also with your client offering to praise your accustomed benevolence, which comes second in no way out of duty. 
For if, as is your custom, you should act kindly, and give to your indebted servant  new pledges of grace, then I will have the accomplishment which to me serves as the greatest joy.
As for the rest, I pay my utmost respects to divine will, for the immense generosity of you, Maecenians and most revered Patrons, for the entire course of happiness, and for all the prosperity which hopefully will endure flourishing for as long as possible. 
Wherever I shall be, I am forever devoted to Your Most Noble Name.

\hfill{Your most devoted follower,}

\hfill{Sven Gust. Sommelius.}

\section{\S I.}

Within Algebra, nothing is more greatly desired than the ability to resolve any equation into its roots. However, although in the pursuit of this matter Mathematical attempts have progressed, a general method of solving quintic equations so far has not been accomplished, unless approximations are sufficient, and this goal remains sought after through the darkest nights. Among different methods, which have been called in reinforcement to remove this difficulty, there is one which is not at all unfavourable, in which a proposed equation is transformed into another, whose solutions are easier to find. 

As an example, assume that the quadratic equation, A, which is to be solved is:
\begin{equation}
  \mbox{A : } z^2+mz+n=0 \ .
\end{equation}
Certainly, nothing is easier than the creation of another Equation B, in which we introduce a new unknown quantity $y$, which satisfies, $z=y-a$.
For if the value of 
$y-a$
is substituted into Equation A, then by this method, Equation A is changed into an equation of this form:
\begin{equation}
  \mbox{B : } y^2-y(-2a+m)+(a^2-ma+n)=0 \ .
\end{equation}
Having made this transformation, I believe it is apparent to anyone that a solution can be found for Equation B with very little difficulty. When I assumed that 
$z=y-a$
as I pleased, we have control over the letter $a$, so nothing else is necessary than to define $a$ in such a way so that the second coefficient becomes zero, and that is when 
$-2a+m=0$
or 
$a=\frac{m}{2}$
Therefore, having substituted in this value for a, Equation B becomes 
\begin{equation}
y^2+\frac{m^2}{4}-\frac{m^2}{2}+n=0 \ 
\mbox{ or }
y^2=\frac{m^2}{4}-n \ .
\end{equation}

In this equation, since we have made it so that the second coefficient disappears, which in equations of this sort is the only intermediary coefficient, by no means will it be hard to find a suitable solution, since it is necessary that:
\begin{equation}
y=\pm \sqrt{\frac{m^2}{4}-n}=\pm \frac{1}{2}\sqrt{m^2-4n}.
\end{equation}
Henceforth, the solution of Equation A, once rearranged, surfaces, since:
$z=y-a=y-\frac{m}{2}$
and so 
\begin{equation}
z=-\frac{m}{2}\pm \frac{1}{2}\sqrt{m^2-4n}. \hspace{2cm} \mbox{Q.E.F. \qed}
\end{equation}
This solution to equation A is the same for all quadratic equations, and is well-known and appears in Algebra books, although in a different way. However, it is clear from this example, not merely what to transform equations means, but also that this can truly resolve equations with great success.

\section{\S II.}

Indeed it is certainly possible to find well-known transformations in all equations with degree of this sort. Truly it is evident, since we only have control over one letter, only one coefficient by this method can be eliminated. Let us assume that our proposed cubic equation is:
\begin{equation}
  \mbox{A: } z^3+mz^2+nz+p=0\
\end{equation}
into which we substitute 
$z=y-a$
and so is transformed into
\begin{equation}
 \mbox{B: } y^3+y^2(-3a+m)+y(3a^2-2ma+n)+(-a^3+ma^2-na+p)=0\ .
\end{equation}
Certainly, in this new equation, the removal of either the second or third coefficients is easy, just as we want: in fact, apart from a very rare case, these coefficients cannot both be reduced at the same time. Because in order to reduce the second coefficient, 
$-3a+m=0, $
that is,
$a=\frac{m}{3}$.
However, to reduce the third coefficient, we need 
$3a^2-2ma+n=0$, or $a=\frac{m}{3}\pm \frac{1}{3}\sqrt{m^2-3n}$.

The solution for $a$, if we were to set
$a=\frac{m}{3}$
is very different from the result if we set $a$ as 
$a=\frac{m}{3}\pm \frac{1}{3}\sqrt{m^2-3n}$, and both results cannot be reached together unless in this case, in which
$\frac{m}{3}=\frac{m}{3}\pm \frac{1}{3}\sqrt[3]{m^2-3n}$, and in which case 
$m^2-3n=0$, or
$n=\frac{m^2}{3}$. 
If therefore we use this equation
$z^3+mz^2+\frac{m^2}{3}z+p=0$, without doubt, taking 
$z=y-\frac{m}{3}$
it is transformed into 
\begin{equation}
y^3+y^2(-m+m)+y(\frac{m^2}{3}-\frac{2m^2}{3}+\frac{m^2}{3})+(-\frac{m^3}{27}+\frac{m^3}{9}-\frac{m^3}{9}+p)=0
\end{equation}
or
$y^3-\frac{m^3}{27}+p=0, $
and in this equation not only the second but also the third coefficient is reduced. However, this case, in which a locus of this kind is possible, is rare to such an extent that when a general method is made to reduce the coefficients, it is not useful for all others that also come under the same category of cubics.

Therefore from the above-mentioned determined values, while one substitution removes more than the other, after one of two has been removed, that term does not remain in the finite equation, in which the other term can be resolved with a solution which is very clearly achievable, and so generally only one term of the transformed equation can be removed. Because when it is as such, a cubic equation can not be resolved up until now through the removal of both intermediate terms. However it cannot be denied that getting rid of the second term provides great ease for solving equations, and that from the written algebra, we reach an adequate solution.

\section{\S III.}
Just as with the quadratic and cubic equations, it is not doubted that in other equations, transformations can be made in such a way. By way of example, let us create this general equation:
\begin{equation}
\mbox{A: } z^\alpha +mz^{\alpha -1}+nz^{\alpha -2}+pz^{\alpha -3}+qz^{\alpha -4}+\ldots=0 \ ,
\end{equation}
in which we assume 
$z=y-a$.
This will become, with the agreement of anyone who is accustomed to the study of limits in mathematics:
\begin{equation}
\begin{split}
\mbox{B: } 
& y^\alpha+y^{\alpha-1}(\frac{-\alpha a}{1}+m)
+y^{\alpha-2} \left( a^2(\frac{\alpha(\alpha-1)}{1\cdot 2})
+ a{\frac{-m(\alpha-1)}{1}}+n \right)+ \\
& + y^{\alpha-3}\left(a^3(\frac{-\alpha(\alpha -1)(\alpha -2)}{1\cdot 2\cdot 3})
 +a^2(\frac{+m(\alpha -1)(\alpha -2)}{1\cdot 2})+a(\frac{-n(\alpha -2)}{1}) \right) + \\
& +y^{\alpha -4}\left(  a^4(\frac{+\alpha (\alpha -1)(\alpha -2)(\alpha -3)}{1\cdot 2\cdot 3\cdot 4}) + 
 a^3(\frac{-m(\alpha -1)(\alpha -2)(\alpha -3)}{1\cdot 2\cdot 3}) \right.+ \\
& \left. \qquad \qquad +a^2(\frac{+n(\alpha -2)(\alpha -3)}{1\cdot 2})+a(\frac{-p(\alpha -3)}{1})+q\right)+ \ldots=0 \ .
\end{split}
\end{equation}
In this equation B, which is more easily observed in this way, for the second term to be removed, 
$-\alpha a+m=0$ 
is substituted in, and also, so that the third term also disappears, 
$a^2(\frac{\alpha (\alpha -1)}{1\cdot 2})-a(m(\alpha -1))+n=0$.
Thus in order to remove the fourth coefficient as required, the substitution is:
\begin{equation}
a^3\left(\frac{-\alpha (\alpha -1)(\alpha -2)}{1\cdot 2\cdot 3}\right)+
a^2\left(\frac{m(\alpha -1)(\alpha -2)}{1\cdot 2}\right)-a(n(\alpha -2))+p=0
\end{equation}
and in general in order to remove the Mth term as required, an equation of degree M needs to be resolved. I believe, with regards to the remaining cases, anyone can easily discover that only one coefficient can be removed, by the same reasoning that has been mentioned already. If this is true, Algebrists finally will come to hope that any sort of equation, freed by the power of transformations, can be resolved by its coefficients, and it must be acknowledged wholly that another method of transformation other than this one has not obtained this to this day, and an expectation of this sort does not mean that an answer will never be found.

\section{\S IV.}
What we have discussed up to this point about the transformation of equations is not new, and generally nothing that Algebra in its youngest form, while still wailing in its cradle, cannot have taught. 
However, no method belonging to the study of Mathematics should have been left behind, untried and extended, from this lovable science by its untiring lovers and worshippers, and this area of Algebra, which deals with transformations of equations, lay unused innocently, doubtless as if forgotten, as Mathematicians will have believed, and until now do believe, that a means by which an equation of the third degree and having two or more coefficients can be solved is certainly impossible using a transformation of this sort.
Indeed either by thinking or otherwise not thinking about this matter, they all wanted to deter the future generations, and they added a demonstration, which is this; They said:
\begin{quote}
Let us create a contradiction; let us assume for example that this cubic equation 
\begin{equation}
\mbox{A: } z^3+mz^2+nz+p=0,\
\end{equation}
can be transformed into another purer cubic equation 
\begin{equation}
\mbox{B: } y^3+a=0,\
\end{equation}
in which both intermediary terms are removed. Therefore the value of $z$ in equation A can be deduced by knowing the value of $y$ in equation B. Since however out of the roots of Equation B only one is real, the other two must be complex; and for Equation A because it holds, due to a lack of loss of generality, such that all cubic equations of this sort adheres to it, and, without doubt, sometimes A is defined by exclusively real roots, which obviously can happen, with the result that a root of equation A, although it is real, nevertheless is determined by a root of equation B, even though it is complex. This says, in effect, that one root, which is real, must be determined by another, which is complex, which creates the contradiction. 
\qed
\end{quote}

We concede with both hands that equation A, which must be kept without loss of generality, can have roots such that they are all real, nor is it denied that equation B can only boast one real root, with the other two always producing the mathematical impossibility. Regarding the possibility of all three roots being real, it is not an easy task. However in the equation B: $y^3+a=0$, one root is always $y+\sqrt[3]{a}=0$ which is real. In order to obtain the other two roots, equation B is divided by the already discovered root, which results in 
\begin{equation}
y^2-y\sqrt[3]{a} +\sqrt[3]{a^2}=0 \ . 
\end{equation}
Thus, out of the remaining two roots, one is
$y-\frac{1}{2}\sqrt[3]{a}+\frac{1}{2}\sqrt[3]{-3\sqrt[3]{a^2}}=0$, and the other is $y-\frac{1}{2}\sqrt[3]{a}-\frac{1}{2}\sqrt[3]{-3\sqrt[3]{a^2}}=0$, both of which are complex. Thus far, in the currently examined demonstration, there has certainly been no lack of mathematical rigour.

When however, it is put forward that it is not possible for a real quantity to be determined by a complex one, this is certainly very general and this phrase \textit{determined by something} is certainly not without suspicion of error. It is ambiguous and it can be defined without doubt such that it can be said non-absurdly that a complex value sometimes can determine a real quantity. If however, this definition of this phrase is taken to be true here, in order to say that a real quantity is refined from parts, of which one or the other can sometimes be complex, then this explanation cannot be doubted, but we come upon in this way a notion that is difficult to understand, which according to the basis of the contradiction cannot be done. Anyone who has delved into the preceding maths or even only common Algebra cannot be unaccustomed to this sort of matter. Let us create, for example, the cubic equation $z^3+mz+p=0$, which by Cardano's formula can be solved into its roots, one of which is this:
\begin{equation}
z-\sqrt[3]{\frac{q}{2}+\sqrt{\frac{q^2}{4}+\frac{p^3}{27}}}-\sqrt[3]{\frac{q}{2}-\sqrt{\frac{q^2}{4}+\frac{p^3}{27}}} \ .
\end{equation}

If it happens that $\frac{q^2}{4}+\frac{p^3}{27}$ is negative, which must happen in some cases, we know $\sqrt{\frac{q^2}{4}+\frac{p^3}{27}}$ becomes complex, and, here, $z$ is the sum of two complex numbers and, the entire argument relies on $z$ being a real quantity. A contradiction of this sort, without doubt, cannot be dealt with in another way than the establishing of the complex numbers, which is $\sqrt[3]{\frac{q}{2}+\sqrt{\frac{q^2}{4}+\frac{p^3}{27}}}$ which is balanced by the other imaginary part which is $\sqrt[3]{\frac{q}{2}-\sqrt{\frac{q^2}{4}+\frac{p^3}{27}}}$ with the result that one cancels out the other. 
There are no exceptions, in which a similar annihilation does not occur, so the value of $z$ in A, although real, nevertheless can be determined from the complex $y$ values in equation B. By the truth of this sought-after demonstration, all contrary assault collapses, such as the demonstration that was supported for a long time more by judgement rather than calculations on the authority of those prestigious mathematicians. Thanks to these ingeniously beautiful efforts to understand the stars in the sky and its sublime flight, never again will this argument be undeservedly overlooked. Certainly if anyone should still remain anxious, we believe they will be placated once they entrust in transformations, which will be used for the abolishment of several terms in any equation.

\section{\S V.}
These results, which were quoted above, apply to all transformations, provided there exists:
\begin{enumerate}
\item A proposed equation A to be solved, which will be transformed.
\item An equation B, which will be created, which contains not only an unknown quantity of equation A, but also a certain new unknown quality, such that the relationship between these two values is defined by the characteristics of equation B. This equation is called the \textit{subsidiary} or \textit{intermediary} equation, since in fact without it a transformation is impossible anyway.
\item A certain algebraic operation, often very difficult, although always possible, by which the unknown quantity which is common to equations A and B is removed. 
\end{enumerate}
Having performed this operation, another equation C appears from the extermination, which tends to be called the \textit{transformed} equation, containing the unknown quantity introduced in equation B, which has the same maximum degree as the other unknown quantity in equation A. In order to further clarify this process, let us repeat it using one of many previously stated examples.
\begin{enumerate}
\item The equation A: $z^3+mz^2+nz+p=0$ is proposed, and is to be transformed;
\item This equation B: $z+a-y=0$ acts as the subsidiary;
\item The letter $z$ is removed from equations A and B, and this transformed equation appears:
\begin{equation}
\mbox{C: } y^3+(-3a+m)y^2+(3a^2-2ma+n)y-a^3+na^2-na+p=0 \ .
\end{equation}
\end{enumerate}
Having considered this, I believe none of the mathematicians can deny that this transformation, described above, and everything else described since the beginning of this dissertation, applies for all cubics. 
We ought to add that we accept, furthermore, to have considered another Subsidiary equation to be applied, which is of a single degree. From this, it happens that we can only control the value of one letter, and thereby only one term from the equation can be removed. Truly, from what I have said of this process, nothing at all impedes a subsidiary equation from having two or higher degree and we will see transformations that happen because of this fact, in which two of more terms can be exterminated from any general equation.

\section{\S VI.}
\begin{enumerate}
\item The equation which is to be transformed is proposed to be:
\begin{equation}
\mbox{A: } z^3+mz^2+nz+p=0,\
\end{equation}
\item This is the subsidiary equation:
\begin{equation}
\mbox{B: } z^2+bz+a+y=0,\
\end{equation}
\item Following the elimination of $z$ from these two equations, this transformed equation appears
\end{enumerate}
\begin{equation}
\begin{split}
\mbox{C: } & y^3+(-mb+m^2-2n+3a)y^2 + \\
&+ (nb^2+(-mn+3p)b+n^2 -2mp-2mba+(2m^2-4n)a+3a^2)y - \\
& -pb^3+mpb^2-npb+p^2 +nb^2a+(-mn+3p)ba+\\
& \qquad \quad + (n^2-2mp)a+(m^2-2n)a^2-mba^2+a^3=0\ .
\end{split}
\end{equation}
In this equation, the two letters $a$ and $b$ are defined however we please, and actually they interact with this equation in such contrasting manners, that they by no means have just the power of a single letter. It is not possible for us to eliminate both the intermediary terms without the use of this trick, which is:
\begin{equation}\begin{split}
1. \quad & \mbox{D: } -mb+m^2-2n+3a=0 \ , \\
2. \quad &\mbox{E: } nb^2+(-mn+3p)b+n^2-2mp-2mba+(2m^2-4n)a+3a^2=0\ .
\end{split}\end{equation}
From equation D, it follows that $a=\frac{mb-m^2+2n}{3}$.
This value of $a$ is substituted into equation E in the place of the same letter $a$, which results in:
\begin{equation}
b^2+\frac{-7mn+9p+2m^3}{3n-m^2}b+\frac{-n^2-6mp-m^4+4m^2n}{3n-m^2}=0 \ .
\end{equation}
The value of this letter $b$ will easily be known, since in order to obtain its value, nothing is needed bar the resolution of a quadratic equation. However, having found out the value of $b$, it is not possible to ignore what $a$ becomes. Furthermore, it follows that, in whatever way letters $a$ and $b$ should be determined, the transformed equation C loses both of its intermediary terms. Q.~E.~F \qed.

Therefore in this way, equation C becomes pure, and is easily resolved. With the known value of $y$ now, the value of $z$ is easily found, after the resolution of the quadratic equation B.
What is very easily understood in this way is that cancelling of one  impossible (imaginary) term by another can happen, as seen in Chapter 4. There are two imaginary roots in equation C, after both intermediary terms disappear, from which the two roots of equation A are derived, even though they are real. In fact, this dependent state cannot be reached unless by using Equation B; however, as long as this equation is simple, that is of a single degree, it does not seem to be something which an imaginary root of equation C can remove.

If however a quadratic equation becomes the intermediary of which both roots can be imaginary, the whole process of explanation becomes very satisfactory. For the imaginary root of the pure equation C can be balanced by the imaginary number, which can be among of the roots of equation B, such that with the complex numbers having been removed one by the other, nothing can stop all the roots of equation A being real and thus generalize for the rest of the cases. 

\section{\S VII.}

From the cubic equations that were transformed to our liking, we are allowed to set forth unobstructed to similarly transform quartic equations by using the quadratic equation. 
Thus this equation is proposed:
\begin{equation}
\mbox{A: } z^4+nz^2+pz+q=0\ ,
\end{equation}
in which the second term is removed, which is necessary in case the matter is pressured by greater difficulties. 
If the second term were to be present, it would certainly not be worth removing it. 
The subsidiary equation is:
\begin{equation}
\mbox{B: } z^2+bz+a+y=0\ .
\end{equation}
Following the removal of the letter $z$, this equation becomes:
\begin{equation}
\begin{split}
\mbox{C: } 
&y^4+(4a-2n)y^3+(6a^2-6na+nb^2+3pb+n^2+2q)y^2 + \\
&+(4a^3-6na^2+(2n^2+4q)a+6pab+2nb^2-npb+2nq+p^2-pb^3-4qb^2)y + \\
& +a^4-2na^3+(n^2+2q)a^2+3pba^2+nb^2a^2-(2nq+2p^3)a - \\
& \qquad \qquad -pnba-4qb^2a-pb^3a+qb^4+nqb^2-qbp+q^2=0\ .
\end{split}
\end{equation}
In this equation, the second and third terms can be removed very easily, if we set:
\begin{enumerate}
\item $4a-2n=0$ or $a=\frac{n}{2}$ and
\item $nb^2+3pb+n^2+2q-\frac{3n^2}{2}=0$.
\end{enumerate}
In truth, the removal of the 2nd and 4th terms seems much more fruitful, since by this, the quartic equation is transformed into a quadratic and is resolved very easily. Regarding the easier solution of the aforementioned matter, we set this substitution which has featured previously, which was how we removed the 2nd and 3rd terms.
Let this be the proposed equation:
\begin{equation}
\mbox{A: } z^4+pz+q=0\ .
\end{equation}
As before, let the subsidiary equation be:
\begin{equation}
\mbox{B: } z^2+bz+a+y=0\ .
\end{equation}
After we substitute out the value of $z$ we will have:
\begin{equation}
\begin{split}
\mbox{C: } y^4+4ay^3+(6a^2+3pb+2q)y^2+(4a^3+4qa+6abp+p^2-pb^3-4qb^2)y\\
+a^4+2qa^2+3pba^2+p^2a-4qb^2a-pba^3+qb^4-qpb+q^2=0\ .
\end{split}
\end{equation}
In this equation if we set $a=0$ then it is necessary that:
\begin{equation}
\mbox{D: } -pb^3-4qb^2+p^2=0\ ,
\end{equation}
in order to remove the 2nd and 4th terms; which is done so that the quartic equation C formally becomes quadratic as desired.
Therefore, the value of $b$ is discovered with the easy resolution of the cubic equation D, and the value of $y$ is then discovered through the completion of the resolution of the quadratic equation C, and moreover, the value of
$z$, which is then itself plucked out of the darkness, cannot be unknown once the quadratic equation B is resolved. \hfill{Q.E.F.}\qed

\section{\S VIII.}

If you were to want to remove the third and fourth terms, we would run into a difficulty that at first glance is neither light nor trivial, and which we ought to treat very carefully; At first we may fear that we will find difficulties in transforming equations of higher degree in many places, and that these will be hindrances to removing the third, fourth and remaining terms. Naturally, in order to simultaneously remove the third and fourth terms of the quartic equation, it is required that:
\begin{enumerate}
\item $6a^2+3pb+2q=0$ or $b=\frac{-6a^2-2q}{3p}$,
\item $4a^3+4qa+6pba-pb^3-4qb^2+p^2=0$,
\end{enumerate}
from which when either $a$ or $b$ are removed, a certain equation is set out, in which the other of these letters is of a higher degree than $z$ or the other unknown values in the proposed equation, thus in order to solve an equation of a lower degree, we need to solve one of higher degree, so bad becomes worse.

Nevertheless, although this difficulty will be present in other equations, it is certain that it can be removed easily in this present occasion. For if at first the second and third terms are removed from the proposed quartic equation, which we have seen to be possible, then a reciprocal transformation is made, and no Mathematician can deny that in this way a quartic equation deprived of its 3rd and 4th terms appears, Q.E.F. \qed 

About the remaining, so that anyone can see very clearly, all equations of any degree can be transformed using an intermediary quadratic equation into another, in which since either the second or the third term is removed, given the possibility of solving a quadratic equation, or the 4th term, given the possibility of the solution of a cubic equation, and thus from now on, in a way I believe nobody can doubt, that in general thanks to this transformation more than two terms can be removed at the same time.

\section{\S IX.}

Therefore, so that we can remove three terms in a certain equation, anyone can see that it is necessary for the intermediary equation to have at least degree 3.
Let the proposed quartic equation be:
\begin{equation}
\mbox{A: } z^4+pz+q=0\
\end{equation}
and the subsidiary equation be:
\begin{equation}
\mbox{B: } z^3+cz^2+bz+a+y=0\
\end{equation}
after the substitution of the letter $z$ we will have:
\begin{equation}
\begin{split}
\mbox{C: } & y^4+(4a-3p)y^3 + (6a^2-9pa+3pbc+4qb+2qc^2+3p^2)y^2 + \\
& +(4a^3-9ba^2+6p^2a+4qc^2a+6pbca+8qba - \\
&	\qquad -pbc^3-4qcb^2-3p^2cb-5pqb+4q^2c+p^2c^3+pqc^2-p^3)y + \\
& +a^4-3pa^3+3p^2a^2+2qc^2a^2+4qba^2+3pbca^2+pqc^2a - \\
&	\qquad -3p^2bca-p^3a+p^2c^3a-4qb^2ca+4q^2ca-5pqba-\\
&	\qquad  - pb^3a+qb^4+3pqcb^2+2q^2b^2-pqbc^3-4q^2c^2b+p^2qb+q^2c^4-pq^2c+q^3=0\ .
\end{split}
\end{equation}
In order to remove all three intermediary terms of this equation C, after duly carried out calculations, none cannot see that it is required that:
\begin{align}
1. \quad &  a=\frac{3p}{4} \ ,\\
2. \quad & \mbox{E:   }  24pbc+16qc^2+32qb-3p^2=0 \ ,\\
3. \quad & \mbox{F:   }  -2pb^3-8qb^2c+2p^2c^3-3p^2bc+4pqc^2-6pqb+8q^2c+p^3=0 \ .
\end{align}
And if we remove $b$ from equations E and F, this equation comes out:
\begin{equation}
\begin{split}
\mbox{G: }&  (24^3p^5-5\cdot 16^3pq^3)c^6+(7\cdot 24^3p^4q-8\cdot 16^3q^4)c^5 +\\
& +540\cdot 24^2p^3q^2c^4+(180\cdot 24p^6+500\cdot 32^2p^2q^3)c^3 +\\
& +(945\cdot 16p^5q+400\cdot 32^2pq^4)c^2+(15\cdot 32\cdot 36q^2p^4+4\cdot 32^3q^5)c+\\
& +7\cdot 32^2p^3q^3-27p^7=0 \ .
\end{split}
\end{equation}
This equation G is of the sixth degree.
It is certainly true that $b$ can be evaluated with the cancellation of the letter $c$ in equations E and F. 
But in fact in this case we do not get an equation of smaller degree. 
It can perhaps be thought that it is not only an equation of sixth degree is present, but actually, under its structure lies an equation of smaller degree, since it is scarcely understood how an expression for the roots of a quartic equation can reveal an expression for roots of a sixth degree equation.

Whatever it is, it seems necessary to solve an equation of the sixth degree in order to remove all the intermediary terms in a quartic equation. 
This is the same difficulty that we encountered earlier in \S.~8, whose solution however is not sufficient as a remedy, and in whose place an extension of this idea will be used. 
We will see however that this difficulty in removing three terms of a certain equation is not always irrefutable.

\section{\S X.}

Let the proposed equation be:
\begin{equation}
\mbox{A: } z^5+pz^2+qz+r=0 \ ,
\end{equation}
in which the second and third terms are to be removed. 
Let the subsidiary equation be:
\begin{equation}
\mbox{B: } z^4+dz^3+cz^2+bz+a+y=0\
\end{equation}
Cancelling out the letter $z$ gives:
\begin{equation}
\begin{split}
\mbox{C: } & y^5+(-3pd-4q+5a)y^4 + \\
&+(3pbc+4qbd+5rb+2qc^2+5rcd-3p^2c+6q^2-4pr+5pqd +
+3p^2d^2-12pda - \\
& \qquad -16qa+10a^2)y^3+\\
& +(-pb^3-4qb^2c-5rb^2d +3p^2b^2+9pbca+12qbda-5rbc^2-3p^2bcd+2pqbc-5pqbd^2+\\
& +15rba+(pr-8q^2)bd-(11rq-3p^3)b+6qc^2a +15rcda+p^2c^3+ pqc^2d+\\
& \qquad + (8rp-4q^2)c^2-9p^2ca+18q^2a+(4q^2-7pr)cd^2-(2qr+3p^3)cd-12pra+\\
& \qquad + 15pqda+(2q^2q+5r^2)c-(p^3+3rq)d^3+9p^2d^2a-(qp^2+5r^2)d^2-\\
& \qquad - (pq^2+rp^2)d-18pda^2+p^4-4q^3+8rpq+10a^3-24qa^2)y^2+ \ldots =0\ .
\end{split}
\end{equation}
In this equation C, the 5th and 6th terms are overlooked, since they are not yet needed.
But in order to remove the second, third and fourth term, it is necessary that:
\begin{align}
1. & \mbox{D: } \quad  a-\frac{3pd+4q}{5}=0; \\
2. & \mbox{E: } \quad 15pbc+20qbd+25rb+10qc^2+25rcd-15p^2c-3p^2d^2-23pqd-2q^2-20rp=0 \ ; \\
\nonumber
3. & \mbox{
The coefficient of the 4th term must also be zero,}\\
\nonumber
& \mbox{ the equation for which we call F.}
\end{align}
If in equation F, we set the value of $a$ equal to $\frac{3pd+4q}{5}$, it is trivial to see that the letters $b$, $c$ and $d$ by this substitution are not raised to a higher degree than before. When however, from equations E and F, either $b$ or $c$ or $d$ are cancelled, this cannot create anything but an equation of sixth degree, which perhaps is not of another form, nor can be changed to a smaller degree by any method. 
Nevertheless, a tentative hope is exhibited of solving this great difficulty.

\section{\S XI.}

If in equation E, being nothing other than the third term of Equation C, instead of $a$ takes the value of $\frac{3pd+4q}{5}$, and we set 
\begin{equation}
b=\alpha d+\zeta,  
\mbox{ and also } 
c=d+\gamma \ . 
\end{equation}
This equation E is turned into this equation:
\begin{equation}
\begin{split}
\mbox{G: } & [(15p+20q)\alpha -3p^2+10q+25r]d^2+\\
& \qquad  [(15p\alpha +20q+25r)\gamma+(15p+20q)\zeta +25r\alpha -15p^2-23pq]d + \\ 
& \qquad +10qy^2 +(15p\zeta -15p^2)\gamma +25r\zeta -2q^2-2rp=0\
\end{split}
\end{equation}
If we set 
\begin{equation}
\alpha=\frac{3p^2-10q-25r}{15p+20q}
\end{equation} 
and 
\begin{equation}
\zeta=\frac{-15p\alpha\gamma-25r\alpha-20q\gamma-25r\gamma+15p^2+23pq}{15p+20q}
\end{equation}
then
\begin{equation}
10q\gamma^2+(-15p^2+15p\zeta)\gamma+25r\zeta-2q^2-2rp=0 \ .
\end{equation}
Anyone can see that in order to find the values of $\zeta$ and $\gamma$ it is not necessary to solve any other equation other than the quadratic one. Therefore, having revealed by this very easy method the values of $\alpha$, $\zeta$ and $\gamma$ and after substituting these in their place in equation G, this equation G is completely reduced, with one term removing another, such that in equation C the third term is also removed. 

Having finished this, we substitute into the 4th term of the same equation C in the place of $a$ the same value of $\frac{3pd+4q}{5}$ and in the place of $b$ the same value of $\alpha d+\zeta$ and also in the place of $c$ the same value $d+\gamma$;
 of course the letters $\alpha$,$\zeta$ and $\gamma$ are determined in the same way that they were be defined before, and we are now looking to remove the additional fourth term from the equation, in which $d$ is the only unknown. 
 However, the maximum degree of $d$ cannot surpass three, so that it is accessible by the resolution of a cubic equation, whose value should then be compatible with the removal of the fourth term in equation C.

 Since therefore, by the definition of $a$ the second term will be removed, and by the definitions of $b$, $c$, $\alpha$, $\zeta$ and $\gamma$ the third term will be removed and by the definition of $d$ the fourth term will be removed, it is true that in this way, the three aforementioned intermediary terms in any equation of 5th degree can be removed. Q.E.F. \qed

In truth, the long extent of this matter and the reason of time forbid all from pursuing this for merit. For the essence of the solution, you can investigate at leisure the nature of the resolution of the cubic and quartic equations, which, by the work of transformations, we saw can become comparable to the rule of Cardano, which concerns greatly these methods.

At present, this will be sufficient and is why we document 
this little piece of work for acceptance.    \hfill{FINEM}

\end{document}